\documentclass[graybox]{svmult}

\usepackage{type1cm}        
\usepackage{makeidx}         
\usepackage{graphicx}        
\usepackage{multicol}        
\usepackage[bottom]{footmisc}

\usepackage{newtxtext}       %
\usepackage[varvw]{newtxmath}       

\makeindex             

\DeclareMathAlphabet{\mathsl}{OT1}{cmss}{m}{sl}
\SetMathAlphabet{\mathsl}{bold}{OT1}{cmss}{bx}{sl}

\newcommand{\ep}{\ensuremath{\varepsilon}}
\newcommand{\om}{\ensuremath{\omega}}

\newcommand{\cL}{\ensuremath{\mathcal L}}
\newcommand{\bbN}{\ensuremath{\mathbb N}}
\newcommand{\bbZ}{\ensuremath{\mathbb Z}}

\DeclareMathOperator{\probv}{\mathbf{P}} 
\DeclareMathOperator{\mean}{\mathbf{E}}

\newcommand{\ldef}{\ensuremath{\mathrel{\mathop:}=}}

\begin{document}

\title*{Heat kernel fluctuations for stochastic processes on fractals and random media}
\titlerunning{Heat kernel fluctuations on fractals and random media}
\author{Sebastian Andres, David Croydon and Takashi Kumagai}
\institute{Sebastian Andres \at Technische Universit\"at Braunschweig, Institut f\"ur Mathematische Stochastik, Universit\"atsplatz 2, 38106 Braunschweig, Germany,
 \email{sebastian.andres@tu-braunschweig.de}
\and David Croydon \at Research Institute for Mathematical Sciences,
Kyoto University, Kyoto 606-8502, Japan, \email{croydon@kurims.kyoto-u.ac.jp}
\and Takashi Kumagai \at Waseda University, 3-4-1 Okubo, Shinjuku-ku, Tokyo 169-8555, Japan, \email{t-kumagai@waseda.jp}}
%
%

%
%
%
%

\maketitle

\vspace{-50pt}
{\it In memory of Professor Robert Strichartz}
\vspace{50pt}

\abstract*{It is well-known that stochastic processes on fractal spaces or in certain random media exhibit anomalous heat kernel behaviour. One manifestation of such irregular behaviour is the presence of fluctuations in the short- or long-time asymptotics of the on-diagonal heat kernel. In this note we review some examples for which such fluctuations are known to occur, including Brownian motion on certain deterministic or random fractals, and simple random walks on various examples of random graph trees, such as the incipient infinite cluster of critical percolation on a regular tree and low-dimensional uniform spanning trees. We also announce some new results that add the one-dimensional Bouchaud trap model to this class of examples.}

\abstract{It is well-known that stochastic processes on fractal spaces or in certain random media exhibit anomalous heat kernel behaviour. One manifestation of such irregular behaviour is the presence of fluctuations in the short- or long-time asymptotics of the on-diagonal heat kernel. In this note we review some examples for which such fluctuations are known to occur, including Brownian motion on certain deterministic or random fractals, and simple random walks on various examples of random graph trees, such as the incipient infinite cluster of critical percolation on a regular tree and low-dimensional uniform spanning trees. We also announce some new results that add the one-dimensional Bouchaud trap model to this class of examples.}

\section{Introduction}

The past 35 years have witnessed an extensive study of heat kernels for stochastic processes on fractals and associated graphs. (See, for instance, \cite{Bar98,Kum14}, and also \cite{Kig01,Str06} for textbooks related to
analysis on fractals.) As a rough illustration of the kind of result that has been proved in this area, one has that if $K$ is a suitably nice non-compact fractal, such as an infinite version of the Sierpinski gasket or Sierpinski carpet, then the corresponding `Brownian motion' $X=\{X_t\}_{t\ge 0}$ on $K$ exhibits sub-Gaussian heat kernel behaviour. More precisely, writing $(p_t(x,y))_{x,y\in K,\:t>0}$ for the transition density of $X$ (typically with respect to a natural choice of Hausdorff measure on $K$), it holds that
\begin{align}
\lefteqn{c_1t^{-\frac{d_f}{d_w}}\exp \left(-c_2 \left(\frac {d(x,y)^{d_w}}{t}\right)^{\frac 1{d_w-1}}\right)\le  p_t(x,y)}\nonumber\\
&\hspace{100pt}\le  c_3t^{-\frac{d_f}{d_w}}\exp \left(-c_4 \left(\frac {d(x,y)^{d_w}}{t}\right)^{\frac 1{d_w-1}}\right),
\label{eq:aron}\end{align}
for all $x,y\in K$, $t>0$, where: $d(\cdot,\cdot)$ is a geodesic distance on $K$, $d_f>0$ is the Hausdorff dimension of $K$ (with respect to $d$), $d_w\ge 2$ is a constant called the walk dimension of the process, and $c_i$, $i=1,\dots,4$, are constants. (Similar estimates are known to hold for simple random walks on graphical versions of such fractals, with the time parameter being restricted to those $t\in 2{\mathbb N}$ satisfying $d(x,y)\le t$.) Crucially, for many fractals, it is known that $d_w>2$, which means the process in question admits anomalous, sub-diffusive space-time scaling. In particular, by integrating \eqref{eq:aron}, one obtains that
\[c_5t^{1/d_w}\le E_x[d(x,X_t)]\le c_6t^{1/d_w},\]
where $E_x$ represents the expectation under the law of $X$ started from $x$, which, when $d_w>2$, is a clear departure from the usual spatial scale of $t^{1/2}$ seen for Brownian motion in Euclidean spaces. (Note that when $d_w=2$, the inequalities at \eqref{eq:aron} represent classical Gaussian estimates.)

For diffusions on domains within Euclidean spaces or on some nice manifolds for which Gaussian estimates hold, more precise asymptotic behaviour of the heat kernel and other objects associated with the spectrum of the diffusions' generators are known, and these are often related to geometric properties of the space. For instance, consider Brownian motion on a bounded domain $\Omega\subset {\mathbb R}^d$, killed upon hitting $\partial \Omega$, and write $\Delta$ for the corresponding Dirichlet Laplacian. Then the well-known result of  Weyl gives that $\rho(x):=\sharp \{\lambda\le x: \lambda ~\mbox{is an eigenvalue of $-\Delta$}\}$, the eigenvalue counting function, satisfies
\begin{equation}\label{weyl}
\lim_{x\to \infty}\frac{\rho(x)}{x^{d/2}}=c_d\,m(\Omega),
\end{equation}
where $m$ is Lebesgue measure on ${\mathbb R}^d$, and $c_d>0$ is a constant depending only upon the dimension of the space. It turns out that on fractals the situation can be quite different, with fluctuations being seen to occur in various places.

It has recently been revealed that fluctuations also occur for stochastic processes on random media, in particular for processes in low dimensions (including certain fractals) and models at criticality. Indeed, in low dimensions, local irregularities in the random medium, in combination with the basic geometry of the space, may affect the process in a non-negligible way that leads to anomalous heat kernel behaviour.
In this paper, we will summarize some results in this area. Specifically, following a discussion of the heat kernel behaviour for two classes of random fractals, we will review heat kernel fluctuations on random trees, namely for random walks on the incipient infinite cluster of a particular critical Galton-Watson tree (i.e.\ the critical percolation cluster on a regular tree conditioned to be infinite) and for random walks on low-dimensional uniform spanning trees. For brevity, we will mainly discuss fluctuations of the on-diagonal parts of the heat kernels on the relevant media. We will also announce some new results concerning fluctuations of the heat kernels in the one-dimensional Bouchaud trap model.

The remainder of the article is organised as follows. In Section \ref{sec2} we review results concerning fluctuations for Brownian motion on deterministic fractals. Sections \ref{sec3} and \ref{sec4} summarize fluctuations for processes on random fractals and random media, respectively.  In Section \ref{sec5}, we announce our recent results concerning the fluctuations of the on-diagonal heat kernel for the one-dimensional Bouchaud trap model.

\section{Fluctuations for Brownian motion on fractals}\label{sec2}

\subsection{Fluctuations of the on-diagonal heat kernel on fractals}

\begin{figure}
\begin{center}
\includegraphics[width=0.6\textwidth]{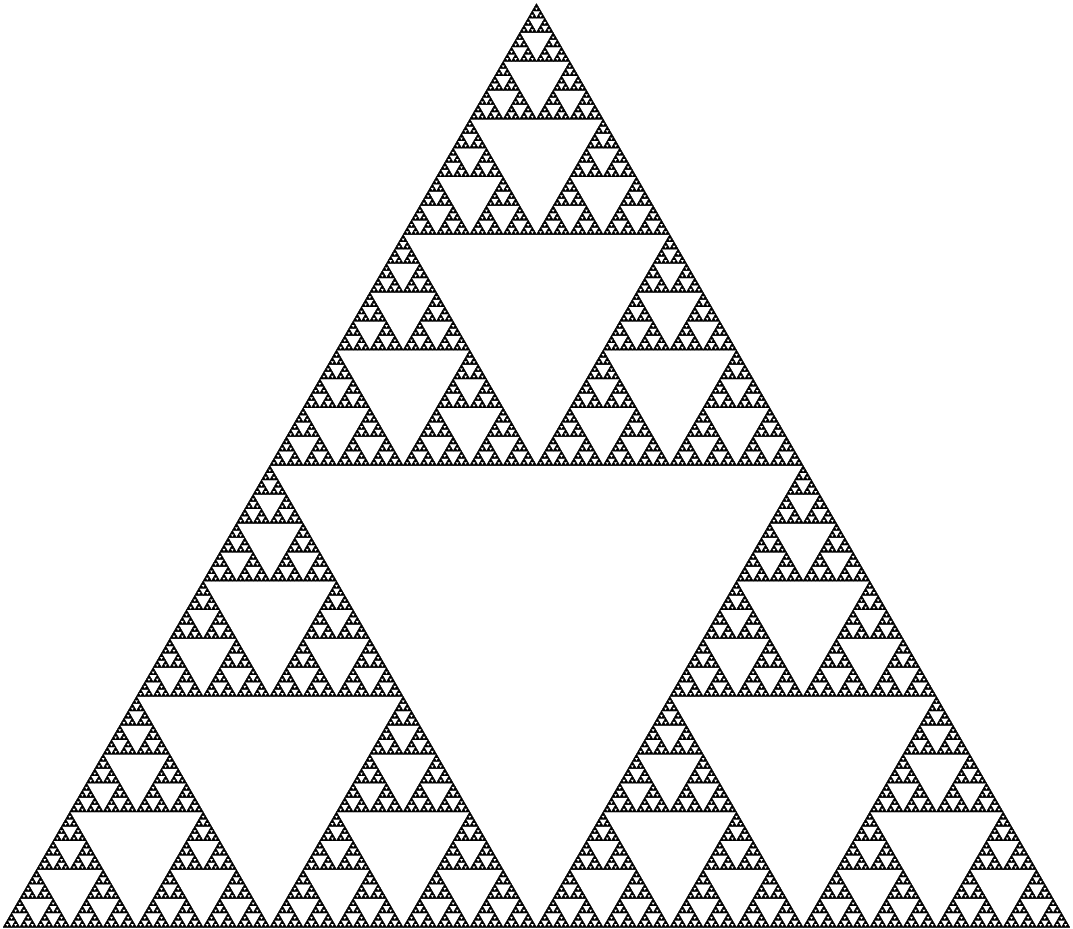}
\end{center}
\caption{The (two-dimensional) Sierpinski gasket.}\label{sgfig}
\end{figure}

Unlike heat kernels on manifolds and regular graphs, heat kernels on fractals exhibit fluctuations. Indeed, in \cite{Kaj13PT}, N.~Kajino shows on-diagonal 
short-time fluctuations for Brownian motion on various fractals. We next give a precise statement in the case of the Sierpinski gasket, as shown in Figure \ref{sgfig}. Note that, in the statement of the result, we write $d_s/2$ in place of the exponent $d_f/d_w$ that appeared in \eqref{eq:aron}; the constant $d_s$ governing the short-time behaviour of the on-diagonal part of the heat kernel is often called the spectral dimension since it also determines the growth of the eigenvalue counting function as in \eqref{weyl} (cf.\ \eqref{eq:disc} below), while in the case of Euclidean domains it matches the actual dimension of the space. For the Sierpinski gasket itself, we have $d_f=\log_2 3$, $d_w=\log_2 5$ and $d_s=2\log_5 3$.

\begin{theorem}{\rm \cite[Theorem 1.1]{Kaj13PT}}
Let $(p_t(x,y))_{x,y\in K,\:t>0}$ be the heat kernel of Brownian motion on the two-dimensional Sierpinski gasket $K$. It is then the case that the limit $\lim_{t\to 0}t^{d_s/2}p_t(x,x)$ does not exist for any $x\in K$.
\end{theorem}

\begin{remark}
(i) In \cite{Kaj13PT}, Kajino obtained the same result for Brownian motion on the $d$-dimensional standard Sierpinski gasket for any $d\geq2$ and on the $N$-polygasket with odd $N\ge 3$ (the latter class including the pentagasket, for example). He also obtained a slightly weaker statement (namely the non-existence of the limit for $\nu$-a.e.\ $x$, for any self-similar measure $\nu$ on the fractal) for Brownian motion on affine nested fractals with at least three boundary points (in the sense of cell boundaries appearing in the construction of the fractal, rather than the topological boundary). See \cite[Theorem 1.2]{Kaj13PT}.\\
(ii) In \cite{Kaj13}, Kajino continued to study the finer behaviour of $t^{d_s/2}p_t(x,x)$ as $t\to 0$. For a certain class of self-similar fractals $K$, he proves that 
$\lim_{t\to 0}t^{d_s/2}p_t(x,x)$
does not exist for generic $x\in K$ if there exists $y,z\in K$ which are not contained in the closure of the boundary of $K$ such that $\limsup_{t\to 0}p_t(y,y)/p_t(z,z)>1$. (Again, boundary should be interpreted in a suitable sense.)\\
(iii) For random walk on the Sierpinski gasket graph, \cite{GW97} established on-diagonal fluctuations for $t\to\infty$. (See also \cite{Woe-00}.)
\end{remark}

\subsection{Other fluctuations}

For Laplace operators and diffusions on fractals, other fluctuations are known, as we will briefly outline next.

\subsubsection*{Spectral properties}

Let $\Delta$ be the Laplace operator on the two-dimensional (compact) Sierpinski gasket $K$. It is known that $-\Delta$ has a compact resolvent, and thus a discrete spectrum. Moreover, defining $\rho$ to be the associated eigenvalue counting function as we did above \eqref{weyl},
\begin{equation}
0<\liminf_{x\to\infty}\frac {\rho (x)}{x^{d_s/2}}<\limsup_{x\to\infty}\frac {\rho(x)}{x^{d_s/2}}<\infty, \label{eq:disc}
\end{equation}
where $d_s=2\log_5 3$ is the spectral dimension mentioned above \cite{FS92}. Notably, unlike \eqref{weyl}, the limit superior and the limit inferior do not coincide in \eqref{eq:disc}. (The periodic nature of the fluctuations are described in \cite{KiLa}.)

Later, in \cite{BaKi}, \eqref{eq:disc} was extended to the Laplace operators on self-similar fractal spaces within quite a general class, including nested fractals. Moreover, it was shown in \cite{BaKi} that the non-existence of the limit is caused by the presence of `many' localized eigenfunctions that produce eigenvalues with high multiplicities. We note that an eigenfunction $u$ of $-\Delta$ is called a localized eigenfunction if $u$ takes the value 0 outside some open subset $O\subset K$ that does not intersect an appropriately-defined fractal boundary. By contrast, as we already noted for Euclidean domains at \eqref{weyl}, the limit at (\ref{eq:disc}) exists with $d_s=d$; no localized eigenfunctions exist in this case.

\subsubsection*{Off-diagonal heat kernels and large deviations}

For Brownian motion $X=\{X_t\}_{t\ge 0}$ on the two-dimensional Sierpinski gasket $K$, fluctuations occur in the Varadhan-type short-time off-diagonal heat kernel asymptotics and in the Schilder-type large deviations principle. Precisely, the following results are known to hold.

\begin{theorem} Let $t_0=2^{1-d_w}=2/5$, and set $\ep_{n,u}=t_0^n u$ for any $u \in [t_0,1)$.\\
(i) {\rm (\cite{Kum97})} It holds that
\begin{equation}
-\lim_{n\to \infty}\ep_{n,u}^{1/(d_w-1)} \log p_{\ep_{n,u}}(x,y)=d(x,y)^{\frac {d_w}{d_w-1}} \, F\left(\frac u {d(x,y)}\right),\qquad \forall x,y\in K,
\label{eq:vard}\end{equation}
where $F$ is a continuous positive non-constant periodic function with period $t_0^{-1}$.\\
(ii) {\rm (\cite{BeKu})} Let $P_x^\ep$ be the law of $B_{\ep \cdot}$ starting at $x$. For each $u\in [t_0,1)$ and Borel measurable $A\subseteq \Omega_x$, where we suppose that $\Omega_x:=\{f\in C([0,T],K):\: f(0)=x\}$ is endowed with the distance induced by the supremum norm, it holds that
\begin{eqnarray}
-\inf_{\phi\in \mbox{\scriptsize Int} (A)} I^u_x(\phi) &\le  &
\liminf_{n\to \infty}\ep_{n,u}^{1/(d_w-1)}\log P_x^{\ep_{n,u}}(A)\nonumber\\
&\le & \limsup_{n\to \infty}\ep_{n,u}^{1/(d_w-1)}\log
P_x^{\ep_{n,u}}(A)
\le  -\inf_{\phi\in \mbox{\scriptsize Cl}( A)} I^u_x(\phi),\label{eq:ldshc}
\end{eqnarray}
where $\mbox{Int} (A)$ and $\mbox{Cl}(A)$ are the interior and closure of $A$, respectively, and $\{I^u_x\}_{u\in [t_0,1)}$ is a collection of rate functions, defined as follows for each $\phi\in \Omega_x$:
\[I^u_x(\phi) \ldef\left\{\begin{array}{ll} \displaystyle{\int_0^T\left(\dot{\phi}(t)\right)^{d_w/(d_w-1)}F\left(\frac {u}{\dot{\phi} (t)}\right) dt },\qquad& \mbox {if $\phi$ is absolutely continuous}, \\ \displaystyle{\infty}, \qquad&\mbox{otherwise},\end{array}\right.\]
with $F$ being the same periodic function as in (i) and $\dot{\phi} (t):=\lim_{s\to t} \frac {d(\phi (s),\phi (t))}{|s-t|}$ for $t\in [0,T]$.
\end{theorem}

For comparison, we note that in the case of Brownian motion on $\mathbb R^n$, \eqref{eq:vard} and \eqref{eq:ldshc} hold independently of $u$, i.e.\ with $F$ being a constant function, and with $d_w=2$. (In fact, there is no need to restrict to a subsequence when considering the Euclidean case.) Hence $I^x_u$ is independent of $u$ and, indeed, one may compute that $I_x(\phi):=\frac 12 \int_0^T\dot{\phi} (t)^2dt$.

In \cite{BeKu}, it was further shown that when $A= \{f\in \Omega_x:\: f(T)=y\}$, then $\inf \{I_x^u(\phi): \:\phi\in A\}$ is attained, independently of $u$, by any path which moves along a geodesic between $x$ and $y$ homogeneously. Thus, regardless of the choice of $u$, geodesics are always `the most probable paths', but the energies (action functionals) of the paths in question depend on time sequences determined by $u$.

\section{Fluctuations of quenched on-diagonal heat kernels on random fractals}\label{sec3}

\begin{figure}
\begin{center}
\includegraphics[width=0.45\textwidth]{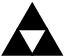}\qquad\includegraphics[width=0.45\textwidth]{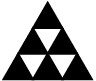}
\end{center}
\caption{The first stage in the iterative construction of SG(2) and SG(3).}\label{SG2SG3}
\end{figure}

Two classes of random Sierpinski gaskets, along with properties of Brownian motions on them, have been considered by B.\ M.\ Hambly, see \cite{Ham92,Ham97}. Both families are based on the choosing randomly from the maps generating the generalised Sierpinski gaskets SG$(\nu)$, $\nu\geq 2$, which are constructed by replacing at each stage a single equilateral triangle with $\nu(\nu+1)/2$ appropriately-arranged smaller ones. See Figure \ref{SG2SG3} for the first stage of the iterative construction of SG(2) and SG(3); note that SG(2) is the standard Sierpinski gasket.

\begin{figure}
\begin{center}
\includegraphics[width=0.45\textwidth]{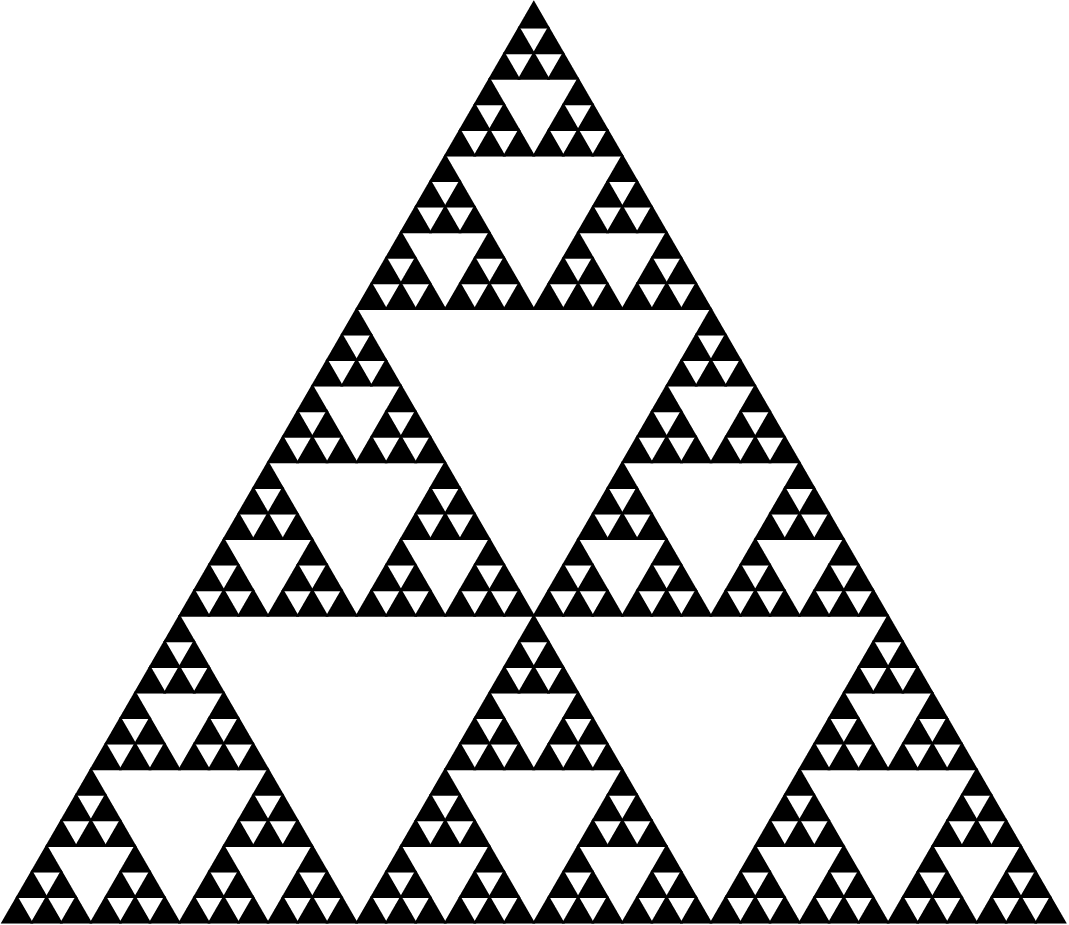}\qquad\includegraphics[width=0.45\textwidth]{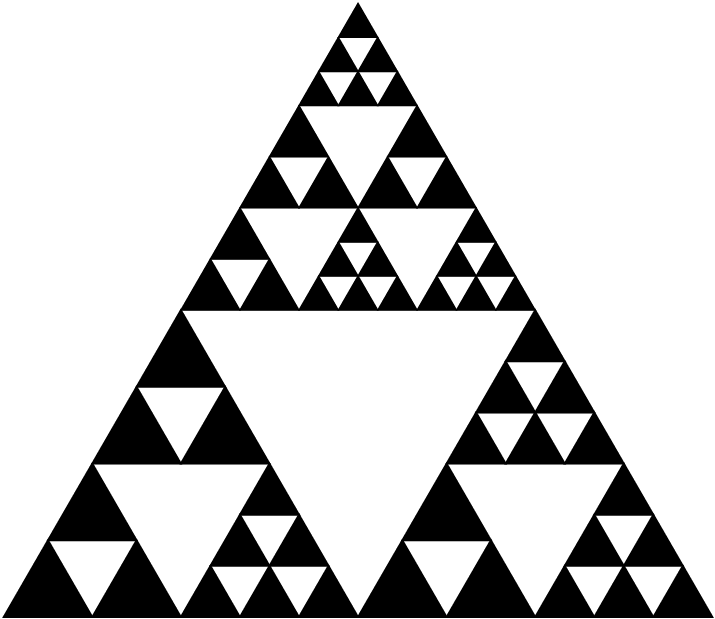}
\end{center}
\caption{Example realisations of a random scale-irregular gasket (left, shown to the fourth level of construction) and a random recursive gasket (right, shown to the third level of construction).}\label{RHGRRG}
\end{figure}

In the model of \cite{Ham92}, one starts with an i.i.d.\ sequence of integers $\{\nu_i\}_{i\in \mathbb N}$, $\nu_i\ge 2$, and then, beginning with a single equilateral triangle, at the $i$-th stage in the construction, divides each remaining equilateral triangle into $\nu_i(\nu_i+1)/2$ sub-triangles, as per the arrangement in the construction of SG($\nu_i$). Via this procedure, one can define a limiting random Sierpinski gasket, which is called the homogeneous random gasket (or the random scale-irregular gasket); see Figure \ref{RHGRRG} for an illustration, and \cite{Ham92} for details. For the Brownian motion on this set, in \cite{Ham92} and \cite[Section 6]{BarHam97}, there are derivations of transition density estimates, which show a sub-Gaussian form as at \eqref{eq:aron} up to some small order error terms, with the leading order exponents being deterministic. See \cite[Theorem 6.1 and Corollary 6.3]{BarHam97} for details, and \cite[Theorem 6.2]{BarHam97} for sufficient conditions on the random sequence $\{\nu_i\}_{i\in \mathbb N}$ (not necessarily restricting to an i.i.d.\ one) for on-diagonal heat kernel fluctuations to occur. Moreover, concerning the associated eigenvalue counting function, \cite[Corollary~7.2]{BarHam97} establishes that this almost-surely satisfies
\[0<\limsup_{x\to\infty}\frac {\rho (x)}{x^{d_s/2}e^{c_1\phi(x)}}<\liminf_{x\to\infty}\frac {\rho(x)e^{c_2\phi(x)}}{x^{d_s/2}}<\infty,\]
where $d_s$ is the appropriately defined spectral 
dimension, as appears
in the on-diagonal part of the heat kernel, $\phi(x)\ldef (\log x \log\log\log x)^{1/2}$, and $c_1$ and $c_2$ are constants. We highlight that $d_s$ is deterministic, meaning that, on a set of probability one, it is independent of $\{\nu_i\}_{i\in \mathbb N}$. In keeping with the content of the present article on heat kernel fluctuations, we further note that the results of \cite{BarHam97,Ham92} are detailed enough to yield that, by replacing the random sequence of integers $\{\nu_i\}_{i\in \mathbb N}$ with a suitably chosen deterministic sequence, one can exhibit examples of scale-irregular gaskets for which the leading polynomial order of the on-diagonal part of the heat kernel is not captured by a single exponent, i.e.\ the spectral dimension does not exist, see \cite{BarHam97} for details.

Another class of random gaskets that is perhaps more closely related to examples of random fractals appearing organically as scaling limits of random graphs (such as those appearing in the next section) are the random recursive gaskets of \cite{Ham97}. In order to understand the construction, it is helpful to first consider the Sierpinski gasket SG$(\nu)$ as (the boundary at infinity of) a $\nu(\nu+1)/2$-ary tree, with vertices in the $n$-th generation of the tree corresponding to the equilateral triangles appearing in the $n$-th level of gasket construction. Generalising this picture, the random recursive gasket is constructed using a random tree corresponding to a simple Galton-Watson tree, whose offspring distribution is supported on integers greater than 2. Vertices in the $n$-th generation of the tree correspond to $n$-th level subsets of the gasket. See Figure \ref{RHGRRG} for an illustration of the initial part of the construction of a random recursive gasket. Notice that, in contrast to the random scale-irregular gasket, the sizes of the different parts of the fractal in the $n$-th stage of the construction can vary.

In \cite{HamKum01}, it was shown that there are fluctuations in the on-diagonal part of the heat kernel and no fluctuations in the leading order of the eigenvalue counting function for random recursive gaskets. For simplicity in presenting these results, consider the case where the random recursive gasket is constructed from only SG$(2)$- and SG$(3)$-type replacements at each stage, i.e.\ the offspring distribution of the associated Galton-Watson branching process is supported on $\{2,3\}$, and let $(\Omega,{\mathcal F},\mathbf{P})$ be the probability space that governs the randomness of the fractal. Furthermore, write $K_\omega$ for a realisation of the random fractal and $\mu_\omega$ for a suitably-defined, statistically self-similar measure on $K_\omega$, where the suffix $\omega$ refers to the underlying element of $\Omega$. In \cite[Theorem 1.1]{HamKum01}, it is shown that there exist $\alpha>0$, $\beta>\beta'>0$, $\beta''>0$ and $c_{1},c_{2},c_3,c_4\in(0,\infty)$ such that the following holds $\mathbf{P}$-a.s.: for $\mu_\omega$-a.e.\ $x\in K_\omega$,
\[c_{1} \leq \limsup_{t\to 0} \frac{t^{\alpha/(\alpha+1)}p^\omega_t(x,x)}
{(\log{|\log{t}|})^{\beta/(\alpha+1)}} \leq c_{2},\]
\[c_3\leq \liminf_{t\to 0}  t^{\alpha/(\alpha+1)}({|\log{t}|})^{\beta''/(\alpha+1)}p^\omega_t(x,x),\quad \liminf_{t\to 0}  \frac{t^{\alpha/(\alpha+1)}p^\omega_t(x,x)}
{(\log{|\log{t}|})^{\beta'/(\alpha+1)}} \leq c_{4}.\]
(The lower bound for the $\liminf$ here was obtained in \cite{Ham97}, and is not expected to be optimal, rather one might conjecture a $\log\log$ term would be sufficient.) Notably, the upper bound for the $\liminf$ and the lower bound for the $\limsup$ establish the on-diagonal heat kernel exhibits log-logarithmic fluctuations, which is a departure from the result for the standard Sierpinski gasket, where bounds of the form \eqref{eq:aron} are seen, and so the fluctuations are only of constant order. We remark that \cite{HamKum01} also considers global fluctuations, i.e.\ the behaviour of $\inf_{x\in K_\omega}p^\omega_t(x,x)$ and $\sup_{x\in K_\omega}p^\omega_t(x,x)$, and here logarithmic corrections are seen. Next, we discuss the eigenvalue counting function of the associated Laplacian (Dirichlet or Neumann). In \cite[Theorem 1.1]{Ham00}, it is shown that there exists a deterministic constant $c>0$ and a strictly positive random variable $W$ with mean one such that
\begin{align*}
 \lim_{x\to\infty} \frac{\rho(x)}{x^{\alpha/(\alpha+1)}} = c W^{1/(\alpha+1)}, \qquad  \text{$\mathbf{P}$-a.s.}
\end{align*}
(In \cite{Ham00}, the results are stated under a more general setting with a certain non-lattice assumption.) In this case, the randomness smooths the object of interest. Indeed, as we described in the preceding section, the eigenvalue counting function of the standard Sierpinski gasket displays periodic fluctuations.

\section{Fluctuations of quenched on-diagonal heat kernels in random media 1: Random Trees}\label{sec4}

The behaviour of the heat kernel on random recursive gaskets is indicative of the on-diagonal heat kernel fluctuations that occur for other kinds of random media, and in particular for various natural classes of random graphs and their scaling limits. (See \cite{Kum14} 
for a survey of random walks on disordered media, including a discussion of anomalous heat kernel estimates.) In this section, we discuss two examples of random media for which heat kernel fluctuations have been proved rigorously.

As before, we will write $\probv$ for the probability measure determining the law of the random medium, and $\mean$ for the corresponding expectation. Given a particular realisation $\omega$ of the random medium, in the context of this section a random (locally finite) graph tree, we write $P^\omega_x$ for the law of the associated (continuous-time) simple random walk $X$ (with unit mean holding times) started from $x$; this is the quenched law of $X$. The quenched transition density is given by
\[p^\omega_t(x,y):=\frac{P^\omega_x(X_t=y)}{\mu^\omega_y},\] 
where $\mu^\om_y$ denotes the number of neighbours of $y$ in the tree associated with the realisation $\om$. (We note that, for a given realisation $\omega$, $\mu^\om_y$ gives the invariant measure of $X$.)

\subsection{Incipient infinite cluster on trees}\label{IIC++}

\begin{figure}
\begin{center}
\includegraphics[width=0.9\textwidth]{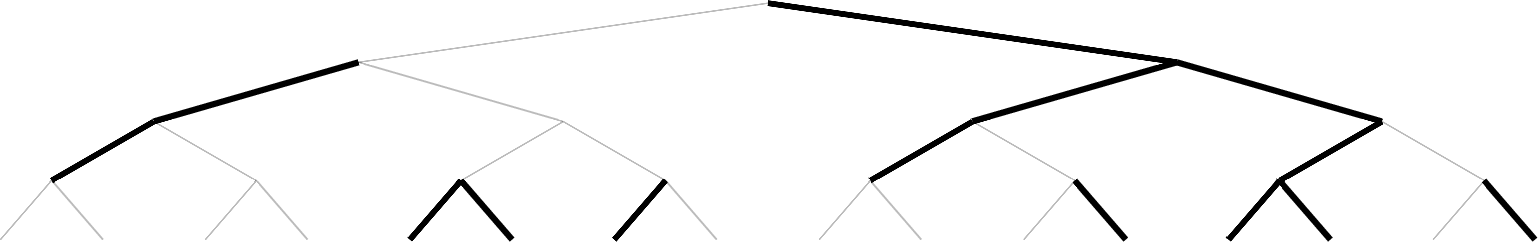}
\end{center}
\caption{Critical percolation on a binary tree. The first elements of the corresponding Galton-Watson branching process are given by $Z_0=Z_1=1$, $Z_2=Z_3=Z_4=2$.}\label{IIC}
\end{figure}

Let $\{Z_n\}_{n\ge 0}$ be a critical Galton-Watson branching process with offspring distribution {\tt Bin}$(n_0,p_c)$, where $p_c=1/n_0$, for some $n_0\in \bbN$. Then the family tree of the process can be understood in terms of a critical percolation process on a regular rooted $n_0$-ary tree ${\mathbb B}$. Indeed, write ${\mathbb B}_N$ for the $N$-th level of ${\mathbb B}$ and ${\mathbb B}_{\le N}$ for the union of the first $N$ levels of ${\mathbb B}$. Then $Z_n=|{\mathcal C}(0) \cap {\mathbb B}_n|$, 
where $0$ is the root of ${\mathbb B}$ and
\[{\mathcal C}(0) :=\{x\in {\mathbb B}: \mbox{ there exists an open path from $0$ to $x$}\},\]
see Figure \ref{IIC}. Since this random tree is almost surely finite, we consider a conditioning procedure that results in a modification of this random tree that extends to infinity. Namely, for $A\subseteq {\mathbb B}_{\le k}$, set
 \[P_0(A):=\lim_{n\to \infty} {\mathbf P}_{p_c}  \left({\mathcal C}(0) \cap {\mathbb B}_{\le k}=A \, | \, Z_n\ne 0\right),\]
where ${\mathbf P}_{p_c}$ is the law of the critical percolation process. It may be checked that the right-hand side is equal to $|A\cap {\mathbb B}_{ k}| \, {\mathbf P}_{p_c} ({\mathcal C}(0) \cap {\mathbb B}_{\le k}=A)$, so the limit exists. Moreover, $P_{0}$ has a unique extension to a probability measure on the set of infinite connected subsets of ${\mathbb B}$ containing $0$. We denote this measure by ${\mathbf P}$, since it will be this probability measure that governs the randomness of the medium in the following. We say that the infinite rooted tree ${\mathcal G}$ chosen according to the distribution ${\mathbf P}$ is called the incipient infinite cluster (IIC). Now, consider a simple random walk on the IIC ${\mathcal G}$. For such a process, the following on-diagonal heat kernel estimates are proved in \cite{BK06}, with $d_f=2$ and $d_w=3$.

\begin{theorem}\label{thm:HKflu}
(i) There exist deterministic constants $c_1,c_2\in(0,\infty)$ and $\alpha>0$ such that,
$\mathbf{P}$-a.s.,
\[c_1t^{-d_f/d_w}(\log\log t)^{-\alpha}\leq p^{\omega}_{t}(0,0)\leq c_2t^{-d_f/d_w}(\log\log t)^{\alpha}\]
for large $t$.\\
(ii) There exists $\beta>0$ such that,  $\mathbf{P}$-a.s.,
\[\liminf_{t\to \infty}  \, (\log\log t)^{\beta} \, t^{d_f/d_w} \, p^{\omega}_{t}(0,0)=0.\]
(iii) There exist constants $c_3,c_4\in(0,\infty)$ such that
\[c_3t^{-d_f/d_w}\leq \mathbf{E}[p^{\omega}_{t}(0,0)]\leq c_4t^{-d_f/d_w}.\]
\end{theorem}

We note that whilst part (i) and (ii) do not establish that heat kernel fluctuations occur $\mathbf{P}$-a.s., they do show that any fluctuations can be at most log-logarithmic in order. Moreover, in combination with the argument used to establish part (iii) they show that log-logarithmic fluctuations occur with positive probability. Heat kernel estimates for critical Galton-Watson trees with more general finite variance offspring distributions are obtained in \cite{FujKum}. Moreover, in the infinite variance case, similar results to the above theorem are derived in \cite{CKtree}, with the distinction that the fluctuations are logarithmic in order in this case.

Finally, we also remark that, for a compact version of the scaling limit of critical Galton-Watson trees in the finite variance case -- the continuum random tree, analogous fluctuation results for the heat kernel were obtained in \cite{CroyCRT}. In fact, the continuum random tree is known to be an example of a random recursive fractal, see \cite{CroyHamCRT}, which goes some way to explaining why the results described above mirror those for the random recursive gaskets of the previous section.

\subsection{Low-dimensional uniform spanning tree}

\begin{figure}
\begin{center}
\includegraphics[width=0.65\textwidth]{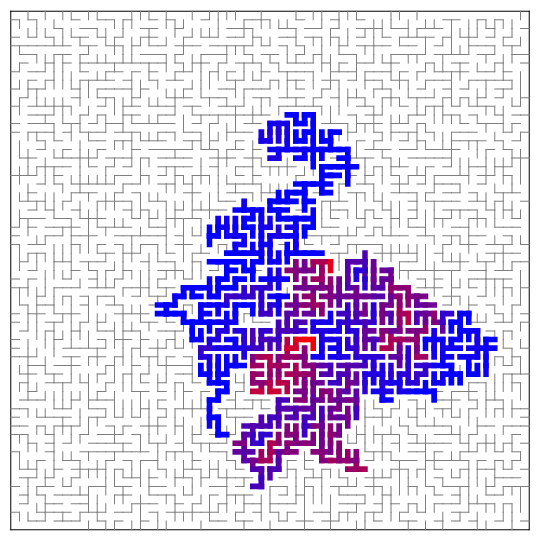}
\end{center}
\caption{The range of a realisation of the simple random walk on uniform spanning tree on a $60\times60$ box (with wired boundary conditions), shown after 50,000 steps. From most to least crossed edges, colours blend from red to blue. Originally presented in \cite{BCK17}, based on code written by Sunil Chhita.}\label{ustfig}
\end{figure}

Let $\Lambda_N:=[-N,N]^2\cap \mathbb{Z}^2$. A loopless connected graph whose vertex set consists of all the elements of $\Lambda_N$ is called a spanning tree of $\Lambda_N$. Let $\mathcal{U}^{(N)}$ be the random graph obtained by picking one of the spanning trees from amongst all possible choices on $\Lambda_N$ uniformly at random, see Figure \ref{ustfig}. The uniform spanning tree (UST) on $\mathbb{Z}^2$, $\mathcal{U}$ say, is the local limit of $\mathcal{U}^{(N)}$ as $N\rightarrow\infty$, see \cite{Pem} for details. Note that it was to describe the potential scaling limit of this model, and the closely connected loop-erased random walk, that O.~Schramm introduced the Schramm-Loewner evolution (SLE) in his celebrated paper \cite{Sch00}.

Now, consider a (discrete-time) simple random walk on $\mathcal{U}$. As for critical Galton-Watson trees, the behaviour of the random walk is anomalous and the quenched on-diagonal heat kernel admits log-logarithmic fluctuations, namely the discrete-time analogue of Theorem~\ref{thm:HKflu} holds with $d_f=8/5$ and $d_w=13/5$. (Note that $8/5=2/(5/4)$; here $2$ is the dimension of the underlying space with respect to the Euclidean metric and $5/4$ is the exponent describing the scaling of the intrinsic metric on $\mathcal{U}$, also with respect to the Euclidean metric.) Parts (i) and (iii) are proved in \cite[Theorems~4.4 and 4.5]{BarMa}, and part (ii) is proved in \cite[Corollary~1.2]{BCK21}, together with the following estimate that establishes the fluctuations occur almost-surely in this case: there exists a deterministic constant $\beta\in(0,\infty)$ such that,  $\mathbf{P}$-a.s.,
\[\limsup_{n\to \infty} (\log\log n)^{-\beta} \, n^{d_f/d_w} \, p^{\omega}_{2n}(0,0)=\infty.\]
Although the estimates and the basic strategies used for proving them are similar to those for the random walk discussed in Subsection~\ref{IIC++}, due to the long-range interaction present in the UST, the proof of Theorem~\ref{thm:HKflu} (in particular (ii)) is much harder. Furthermore, we note that we state here the result for the discrete-time walk since that is what is considered in \cite{BCK21,BarMa}, but it would be straightforward to deduce the corresponding continuous-time result from the same essential estimates on the geometry of the space.

\begin{remark} Recently this problem has also been studied for random walks on the UST on $\mathbb{Z}^3$. With suitably modified exponents, parts (i) and (iii) of Theorem~\ref{thm:HKflu} are proved in \cite[Theorem 1.9]{ACHS21},  and the fluctuation result of Theorem~\ref{thm:HKflu}(ii) is proved in
\cite[Theorem 1.2]{ShirW}.
\end{remark}

\section{Fluctuations of quenched on-diagonal heat kernels in random media 2: Bouchaud trap models}\label{sec5}

In this section, we announce some heat kernel fluctuation results for another kind of random media, namely the one-dimensional Bouchaud trap model, which are obtained in the ongoing work \cite{ACK}.

To introduce the model, let $\tau=\{\tau_x\}_{x\in \bbZ}$ be i.i.d.\ random variables with distribution given by
\[\probv\left(\tau_x >u\right) = u^{-\alpha}\]
for $u\geq 1$, where $\alpha>0$ is some fixed constant; this will represent the random environment in the current setting. Given a realisation of $\tau$, let $X=(X_t)_{t\geq 0}$ be the continuous-time Markov chain with generator
\[\cL f(x) = \frac{1}{2\tau_x} \sum_{y\sim x} \left(f(y)-f(x) \right);\]
this is the one-dimensional (symmetric) Bouchaud trap model. The dynamics of $X$ are as follows: $X$ waits at a vertex $x$ an exponentially-distributed time with mean $\tau_x$, and, at this time, it jumps to one of the two neighbours of $x$, chosen uniformly at random. It is easy to see that the measure $\mathbf{\tau}\ldef \sum_{x\in \bbZ} \tau_x \delta_x$ is invariant for the reversible Markov chain $X$. As before, we will write $P_x^\omega$ for the quenched law of $X$ started from $x$, where again we suppose $\omega$ is a variable that determines the environment, so in particular we have that $\tau=\tau(\omega)$. Moreover, for this model, we will also consider
\[\mathbb{P}(\cdot)=\int P_0^\omega(\cdot) \, \mathbf{P}(d\omega),\]
which is typically called the annealed law of $X$ started from 0.

We first recall the known scaling limits for the process $X$. The topology considered in the following theorem is the usual Skorohod $J_1$ topology on the space of c\`{a}dl\`{a}g paths $D([0,\infty),\mathbb{R})$, and we omit to mention a (non-trivial) deterministic constant in the time-scaling of the limiting process.

\begin{theorem}
(i) If $\alpha\in(0,1)$, then $(\varepsilon X_{t/\varepsilon^{1+1/\alpha}})_{t\geq 0}$ converges in distribution under the annealed law $\mathbb{P}$ to the FIN diffusion (for Fontes-Isopi-Newman, see \cite[Definition~1.2]{FIN}) with parameter $\alpha$.\\
(ii) If $\alpha=1$, then $(\varepsilon X_{t\log(\varepsilon^{-1})/\varepsilon^{2}})_{t\geq 0}$ converges in distribution under the annealed law $\mathbb{P}$ to Brownian motion.\\
(iii) If $\alpha>1$, then, for $\mathbf{P}$-a.e.\ realisation of the environment $\tau$, $(\varepsilon X_{t/\varepsilon^{2}})_{t\geq 0}$ converges in distribution under the quenched law $P^\omega_0$ to Brownian motion.
\end{theorem}

We note that part (i) is contained within \cite[Theorem~3.2(ii)]{BCCR}. We also refer to \cite{FIN} for the introduction of the limiting process and a scaling limit for a single marginal of the process, and \cite{BC05} for a related result in the non-symmetric case. Part (iii) is contained within  \cite[Theorem~3.2(ii)]{BCCR}. Part (ii) is described in \cite[Remark~3.3]{BCCR}, and can also be checked by applying the resistance scaling techniques of \cite{Croyres,CHKres}. Indeed, via such an approach a similar result (i.e.\ convergence to Brownian motion under the annealed law, with logarithmic terms in the scaling) has been obtained in the more challenging case of the random walk on the range of four-dimensional random walks in \cite{Cshir}.

To capture the leading order term in the on-diagonal part of the heat kernel of $X$, we define
\[\phi_\alpha(t):=\left\{
                    \begin{array}{ll}
                      t^{-\frac{1}{1+\alpha}}, & \hbox{if $\alpha\in(0,1)$;} \\
                      t^{-\frac{1}{2}}(\log t)^{-\frac{1}{2}}, & \hbox{if $\alpha=1$;} \\
                      t^{-\frac12}, & \hbox{if $\alpha>1$.}
                    \end{array}
                  \right.\]
The following theorem, which shows quenched fluctuations of the heat kernel about $\phi_\alpha$ when $\alpha\leq 1$, is due to \cite{ACK}. We note that for the short-time heat kernel asymptotics of the limiting FIN diffusion, similar results to those contained in part (i) of the result were obtained in \cite{CHK}.

\begin{theorem}\label{thm:QHK}
(i) If $\alpha\in(0,1)$, then it $\probv$-a.s.\ holds that
\[\limsup_{t\rightarrow\infty}\frac{p_t^\omega(0,0)}{\phi_\alpha(t)(\log\log t)^{\frac{1-\alpha}{1+\alpha}}}\in(0,\infty),\]
and
\[\liminf_{t\rightarrow\infty}\frac{p_t^\omega(0,0) \, (\log t)^{\theta}}{\phi_\alpha(t)}=\begin{cases}
0,  & \text{if $\theta=\frac{1}{1+\alpha}$},\\
\infty, & \text{if $\theta>\frac{3}{\alpha}$.}
\end{cases}\]
(ii) If $\alpha=1$, then it $\probv$-a.s.\ holds that
\[\limsup_{t\rightarrow\infty}\frac{p_t^\omega(0,0)}{\phi_\alpha(t)}\in(0,\infty),\]
and
\[\liminf_{t\rightarrow\infty}\frac{p_t^\omega(0,0) \, (\log\log t)^{\theta}}{\phi_\alpha(t)}=\begin{cases}
0,  & \text{if $\theta=\frac12$},\\
\infty, & \text{if $\theta>3$.}
\end{cases}\]
(iii) If $\alpha>1$, then there exists a constant $\sigma^2>0$ such that, $\probv$-a.s., for every $x_0\in[0,\infty)$ and compact interval $I\subset (0,\infty)$,
    \[\lim_{\lambda\rightarrow\infty}\sup_{|x|\leq x_0}\sup_{t\in I}\left|\lambda p^\omega_{\lambda^2t}\left(0,\lfloor \lambda x\rfloor\right)-\frac{1}{\sqrt{2\pi\sigma^2t}}e^{-\frac{x^2}{2\sigma^2 t}} \right|=0.\]
    In particular,
    \[\lim_{t\rightarrow\infty}\phi_\alpha(t)^{-1}p_t^\omega(0,0)=\frac{1}{\sqrt{2\pi\sigma^2}}.\]
\end{theorem}

Let us very briefly mention the idea of the proof, which is based on standard ideas in the area that relate heat kernel estimates to estimates on the resistance and volume growth of a space. (A framework for deriving such estimates that is broad enough to include random recursive gaskets, the continuum random tree and the FIN diffusion is presented in \cite{CroyLMS}.) In the setting of the one-dimensional Bouchaud trap model, the resistance metric is easy to understand, since it simply coincides with the Euclidean metric. As for the volume growth of the invariant measure, consider
\begin{align*}
V(x,n) \ldef \sum_{y=x-n}^{x+n} \tau_y, \qquad x\in \bbZ, \, n\in \bbN.
\end{align*}
In particular, as a sum of i.i.d.\ random variables, the asymptotic behaviour of $V(0,n)$ can be understood in terms of classical results of probability theory. When $\alpha>1$, so the random variables in the sequence $\tau$ have a mean, the law of large numbers immediately implies that, $\mathbf{P}$-a.s., $V(0,n)\sim cn$, for some deterministic constant $n$. Indeed, with a little more work, one can check that a suitably-scaled version of the invariant measure converges to Lebesgue measure on the line. Since there are no fluctuations of the volume in this case, there are no fluctuations of the heat kernel, as seen in part (iii) of the above theorem. On the other hand, when $\alpha<1$, the random variables in the sequence $\tau$ fall into the domain of attraction of an $\alpha$-stable law, and so we see subordinators in the distributional scaling limit of $V(0,n)$. Again, with more work, one can describe the distributional scaling limit of the invariant measure, this time in terms of a Poisson random measure. Now, in such a regime, there is no quenched convergence of the volume, but rather we see $\mathbf{P}$-a.s.\ fluctuations about the leading order scaling factor. It is this that leads to the fluctuations of part (i) of Theorem \ref{thm:QHK}. Note that, in this case, we see asymmetric fluctuations, with upper fluctuations of log-logarithmic order, but lower ones of logarithmic order. As for $\alpha=1$, this falls somewhere between the two other cases, with a deterministic distributional limit for the volume, but quenched fluctuations about this. Thus, apart from the intermediate case of $\alpha=1$, this example captures the basic heuristic that one might expect to see heat kernel fluctuations when the scaling limit of the random media is itself random.

Finally, we note that, as for the volume measure, the heat kernel of $X$ shows smoother behaviour in terms of its distributional scaling \cite{ACK}; this is reminiscent of the result for critical Galton-Watson trees that was presented above as Theorem~\ref{thm:HKflu}(iii).

\begin{theorem}
(i) For every $\alpha>0$,
\[\lim_{\lambda\rightarrow\infty}\liminf_{t\rightarrow\infty}
\probv\left(\lambda^{-1}\leq \phi_\alpha(t)^{-1}p_t^\omega(0,0)\leq \lambda\right)=1.\]
(ii)  For every $\alpha>0$, there exists an $\varepsilon>0$ such that
\[0<\liminf_{t\rightarrow\infty} \phi_\alpha(t)^{-\varepsilon}\mean\left(p_t^\omega(0,0)^\varepsilon\right)
\leq \limsup_{t\rightarrow\infty}  \phi_\alpha(t)^{-\varepsilon}\mean\left(p_t^\omega(0,0)^\varepsilon\right)<\infty.\]
Moreover, if $\alpha>\frac{3}{2}$, then we may take $\varepsilon=1$.
\end{theorem}

The restriction on the moments in part (ii) above relates to the integrability of certain powers of the volume that appear in the heat kernel estimates. Similar bounds were obtained in \cite{CKtree} for random walks on critical Galton-Watson trees with infinite variance offspring distributions.

\begin{acknowledgement}
This research was supported by JSPS Grant-in-Aid for Scientific Research (A) 22H00099, JSPS Grant-in-Aid for Scientific Research (C) 19K03540, and the Research Institute for Mathematical Sciences, an International Joint Usage/Research Center located in Kyoto University.
\end{acknowledgement}

\bibliographystyle{abbrv}
\bibliography{lit-Bob}
\end{document}